\documentclass[12pt,twoside]{article}
\usepackage{amssymb,amsmath,amsthm,latexsym,amscd}
\usepackage{amscd,graphics,graphicx,cite}
\usepackage[cp1251]{inputenc}
\usepackage[english]{babel}

\newtheorem{theorem}{Theorem}
\newtheorem{lemma}{Lemma}

{\theoremstyle{definition}

\newtheorem{remark}{Remark}

}

\DeclareMathOperator{\Rep}{Rep}

\DeclareMathOperator{\Imm}{Im}
\DeclareMathOperator{\diag}{diag}
\DeclareMathOperator{\End}{End}
\DeclareMathOperator{\Rad}{Rad}

\newcommand{\src}[1]{\stackrel{\circ}{#1}}
\newcommand{\srb}[1]{\stackrel{\bullet}{#1}}
\newcommand{\seq}[2]{#1_1,#1_2,\ldots ,#1_#2}

\newcommand{\wt}[1]{\widetilde{#1}}
\newcommand{\wh}[1]{\widehat{#1}}
\newcommand{\mc}[1]{\mathcal{#1}}
\newcommand{\mb}[1]{\mathbb{#1}}
\newcommand{\mf}[1]{\mathfrak{#1}}
\newcommand{\ov}[1]{\overline{#1}}
\newcommand{\oa}[1]{\overrightarrow{\mkern -3mu #1\mkern 3mu}}
\newcommand{\RepQH}{\Rep(Q,\mathcal{H})}
\newcommand{\ReposQH}{\Rep_{os}(Q,\mathcal{H})}
\newcommand{\ReposQnH}{\Rep_{os}(Q_n,\mathcal{H})}
\newcommand{\mcl}[1]{\multicolumn{1}{|c}{#1}}
\newcommand{\mcr}{\multicolumn{1}{|l}{}}

\begin{document}

\begin{center}
\Large \bf
On the indecomposability in the category of representations of a
quiver and in its subcategory of orthoscalar representations
\end{center}

\begin{center}
\large\bfseries S.~A.~Kruglyak~$^\dag$, L.~A.~Nazarova, \fbox{A.~V.~Roiter}
\end{center}

\noindent
$^\dag$~{\small E-mail: red@imath.kiev.ua}\bigskip

Representations of quivers of the finite and tame types are classified up to
equivalence in the papers \cite{Gab1,Naz1}.

It is naturally to classify representations of quivers in the category of Hilbert
spaces up to the unitary equivalece \cite{Roi1,KruRoi1}.
One can regard the category of such represenations as a subcategory in the
category of all representations, and at that objects, which are indecomposable
in the subcategory, become in general decomposable in the ``larger'' category.
It does not happen with indecomposable locally scalar representations
\cite{KruRoi1} of a quiver.

A proof of this fact and its applications are in present paper.

\bigskip

{\bf 1.} A quiver $Q$ with the set of vertices $Q_v,\; |Q_v|=N$ and
the set of arrows $Q_a$ is called \emph{separated} if $Q_v=\src{Q}\sqcup\srb{Q}$,
and for any $\alpha\in Q_a$ its tail $t_\alpha\in\src{Q}$
and head $h_\alpha\in\srb{Q}$. A quiver $Q$ is \emph{single} if
all its arrows are single
(i.~e. if $\alpha\neq\beta$ then either $t_\alpha\neq t_\beta$ or
$h_\alpha\neq h_\beta$). Vertices from $\src{Q}$ are called
\emph{even,} and from $\srb{Q}$ are \emph{odd.}

Let $m=|\srb{Q}|,\; n=|\src{Q}|,\; \srb{Q}=\{\seq{i}{m}\},\;
\src{Q}=\{\seq{j}{n}\}$. A representation $T$ of a quiver $Q$
associates a finite-dimensional linear space $T(i)$ to a vertex
$i\in Q_v$, and a linear map $T_{ij}: T(j)\to T(i)$ to an arrow
$\alpha : j\to i,\; \alpha\in Q_a$.

Representation $T$ of single separated quiver for fixed
bases of spaces $T(i),\; i\in G_v$ can be associated with a matrix,
separated by $m$ horizontal and $n$ vertical bars, i.~e. the matrix
$$
T=\begin{bmatrix}
T_{i_l,\, j_k}
\end{bmatrix}
\substack{k=\ov{1,n}\\l=\ov{1,m}}.
$$
Here we assume that $T_{i_l,\, j_k}=0$, if there does not exist any such
$\alpha\in Q_a$ that $t_\alpha = j_k,\; h_\alpha = i_l$.

Let $\Rep Q$ be the category of representations of a quiver $Q$, which objects
are representations, and morphism of a representation $T$ to a representation
$\wt{T}$ is defined as a family of linear maps $C=\{C_i\}_{i\in Q_v}$
such that for each
$\alpha\in Q_v$ with $t_\alpha = j,\; h_\alpha = i$ the diagram
\begin{equation}\label{eq01}
\begin{CD}
T(j) @>T(\alpha)>> T(i)\\
@V C_j VV @V C_i VV\\
\wt{T}(j) @>\wt{T}(\alpha)>> \wt{T}(i)
\end{CD}
\end{equation}
\noindent is commutative, i.~e. $C_iT_{ij}=\wt{T}_{ij}C_j$.

Let representations $T,\,\wt{T}$ of a separate single quiver
are defined in the matrix form by the matrices
$$
T=\begin{bmatrix}
T_{i_l,\, j_k}
\end{bmatrix}
\substack{k=\ov{1,n}\\ l=\ov{1,m}}
\quad\text{and}\quad
\wt{T}=\begin{bmatrix}
\wt{T}_{i_l,\, j_k}
\end{bmatrix}
\substack{k=\ov{1,n}\\ l=\ov{1,m}}
$$
Let $C: T\to \wt{T}$ be a morphism of representations, $C=\{C_i\}_{i\in Q_v}$.
Let us introduce the matrices
$A=\diag\{C_{i_1},C_{i_2},\ldots,C_{i_m}\},\;
B=\diag\{C_{j_1},C_{j_2},\ldots,C_{j_n}\}$.
Then commutativity of the diagrams \eqref{eq01} implies
\begin{equation}\label{eq02}
AT=\wt{T}B
\end{equation}
We will further say that $C=(A,B)$.

Let $\mc{H}$ be the category of unitary (finite-dimensional Hilbert) spaces.
Denote as $\RepQH$ the subcategory in $\Rep Q$ which objects are
representations $T$ for which $T(i)$ are unitary spaces ($i\in Q_v$) and
morphisms $C: T\to \wt{T}$ are those of morphisms in $\Rep Q$ for which,
except \eqref{eq01}, the following diagrams are also commutative
\begin{equation}\label{eq03}
\begin{CD}
T(j) @<T(\alpha)^*<< T(i)\\
@V C_j VV @V C_i VV\\
\wt{T}(j) @<\wt{T}(\alpha)^*<< \wt{T}(i)
\end{CD}
\end{equation}
\noindent i.~e. the following equality holds:
\begin{equation}\label{eq04}
BT^*=\wt{T}^*A
\end{equation}
Two representations $T,\,\wt{T}$ from $\Rep Q$ are equivalent in $\Rep Q$, if
there exists an invertible morphism $C: T\to \wt{T}$ in $\Rep Q$.

Two representations $T,\, \wt{T}$ from $\RepQH$ are equivalent $\RepQH$, if
there exists an invertible morphism $C: T\to \wt{T}$ in $\RepQH$.
It can be shown that $T,\,\wt{T}$ are equivalent in $\RepQH$ if and only if
they are unitary equivalent (see, for instance, \cite{Roi1}), i.~e. the
invertible morphism $C$ can be chosen as consisting of unitary matrices $C_i$.

Denote
\begin{gather*}
\oa{T}_i=\;
\begin{matrix}
\cline{1-4}
\mcl{T_{i,\,j_1}} & \mcl{T_{i,\,j_1}} & \mcl{\cdots} &
\mcl{T_{i,\,j_n}} & \mcr \\
\cline{1-4}
\end{matrix} \\\\
T_i^\downarrow =\;
\begin{matrix}
\cline{1-1}
\mcl{T_{i_1,\,j}} & \mcr \\
\cline{1-1}
\mcl{T_{i_2,\,j}} & \mcr \\
\cline{1-1}
\mcl{\vdots} & \mcr \\
\cline{1-1}
\mcl{T_{i_m,\,j}} & \mcr \\
\cline{1-1}
\end{matrix} \\
\oa{T}_i : \sum_{k=1}^n \oplus T(j_k)\to T(i),\\
T_j^\downarrow : T(j) \to \sum_{l=1}^m \oplus T(i_l).
\end{gather*}
Representation $T$ of a separate single quiver $Q$ from the category
$\RepQH$ is called \emph{orthoscalar} (in~\cite{KruRoi1} such representations
are called locally scalar) if each $i\in Q_v$ is associated with a real
positive number $\chi_i$, and the following conditions hold:
\begin{gather}
\oa{T}_i\cdot\oa{T}_i^* = \chi_iI_i\;\text{for}\;i\in\srb{Q},\nonumber\\
T_j^{\downarrow *}\cdot T_j^\downarrow = \chi_jI_j\;\text{for}\;j\in\src{Q},
\label{eq05}
\end{gather}
\noindent here $I_i$ is the matrix of identity operator in $T(i)$.

Define the category $\ReposQH$ as a full subcategory in $\RepQH$,
which objects are orthoscalar representations of a quiver $Q$.

Let us associate two $N$-dimensional vectors ($N=m+n$)
to an orthoscalar representation $T$ of a separate single quiver $Q$:
a \emph{dimension of representation $T$}:
$d=\{d(j)\}_{j\in Q_v}$, where $d(j)=\dim T(j)$, and
a \emph{character of representation $T$}: $\chi = \{\chi(j)\}_{j\in Q_v}$,
$\chi(j)=\chi_j$ are defined above in \eqref{eq05}.
It is easy to see that
$$
\sum_{l=1}^m d(i_l)\chi(i_l)=
\sum_{k=1}^n d(j_k)\chi(j_k)=
$$
\noindent (the sum of squares of the lengths of rows in a matrix of a
representation $T$ equals to the sum of squares of the lengths of columns).

Two orthoscalar representations $T$ and $\wt{T}$ are equivalent if there exist
such unitary matrices
$U=\diag\{U_{i_1},U_{i_2},\ldots,U_{i_m}\}$ and
$V=\diag\{V_{j_1},V_{j_2},\ldots,V_{j_n}\}$ that
\begin{gather}\label{eq06}
UT=\wt{T}V,\;\text{or}\nonumber\\
\wt{T}_{i_l,\,j_k}=U_{i_l}T_{i_l,\,j_k}V_{j_k}^*.
\end{gather}
Here the matrices $U_{i_l}$ have dimension $d_{i_l}\times d_{i_l}$, and
$V_{j_k}$~--- dimension $d_{j_k}\times d_{j_k}$.

\begin{remark}\label{re01}
Note that for an arbitrary quiver $Q$ we can construct $*$-quiver $\wh{Q}$,
adding to each arrow $\alpha_{ij}: j\to i$ an additional arrow
$\alpha_{ij}^*: i\to j$. At that the category $\ReposQH$
and the category of orthoscalar $*$-representations of a $*$-quiver $\wh{Q}$
are identified.
\end{remark}

\begin{remark}\label{re02}
Let $C=(A,B)$ be a morphism of a representation $T$ to a representation $\wt{T}$
in the category $\Rep Q$, i.~e. the equality \eqref{eq02} holds:
$$
AT=\wt{T}B
$$
\noindent and $A,\,B$ are unitary operators, then $C=(A,B)$ is also a
morphism of representation $T$ to representation $\wt{T}$
in the category $\RepQH$, i.~e. the equality \eqref{eq04} holds:
$$
BT^*=\wt{T}^*A
$$
Indeed, \eqref{eq02} implies $T^*A^*=B^*\wt{T}^*$ or,
considering the unitarity of $A$ and $B$, we have
$T^*A^{-1}=B^{-1}\wt{T}^*$. Therefore $BT^*=\wt{T}^*A$.
\end{remark}

\begin{remark}\label{re03}
Let $C=(A,B)$ be a morphism of a representation $T$ to representation $T$
(endomorphism of representation $T$) in the category $\Rep Q$, i.~e.
\begin{equation}\label{eq07}
AT=TB
\end{equation}
\noindent and $A,\,B$ are self-adjoint operators, then $C=(A,B)$ is also an
endomorphism of representation $T$ in the category $\RepQH$, i.~e.
\begin{equation}\label{eq08}
BT^*=T^*A
\end{equation}
Indeed, \eqref{eq08} is obtained from \eqref{eq07} by adjunction.
\end{remark}

\bigskip

{\bf 2.} Representation $T$ of a quiver $Q$ is \emph{faithful} if $T(i)\neq 0$
for all $i\in Q_v$ ($d(i)\neq 0$). \emph{Support} of representation $T$
is the set $Q_v^T=\{i\in Q_v \,|\, T(i)\neq 0\}$.
The character of orthoscalar representation is defined uniquely on the support
$Q_v^T$ (and not uniquely out of it). If a representation is faithful then
the character of orthoscalar representation is uniquely defined and denoted as
$\chi_T$; in general we denote as $\{\chi_T\}$ the set of all characters of
representation $T$. Evidently, if $T$ and $\wt{T}$ are equivalent, then
$\{\chi_T\}=\{\chi_{\wt{T}}\}$.

Decomposable representations are defined in the natural way (at that, if
$T=T_1\oplus T_2$ in the category $\RepQH$, then $T_1(i)\oplus T_2(i)$ means
orthogonal sum of unitary spaces), and if $T,\,T_1,\,T_2$ are
faithful orthoscalar representations, then $\chi_{T_1}=\chi_{T_2}=\chi_T$.

Representation $T$ is called \emph{Schur representation}
in the category $\Rep Q$
(resp. $\RepQH,\, \ReposQH$), if its ring of endomorphisms
is one-dimensional (resp. is isomorphic to $\mb{C}$).

Obviously, if $T$ is Schur representation then $T$ is indecomposable
(in the respective category).

\begin{remark}\label{re04}
If $T$ is indecomposable representation in the category $\mc{R}=\RepQH$ then
$T$ is Schur representation.

Indeed, the algebra $\mf{A}=\End_{\mc{R}}T$ is finite-dimensional $*$-algebra.
If $C=(A,B)\in\End_{\mc{R}}T$ then $C^*=(A^*,B^*)$.
If $C\in\Rad\mf{A}$ then $CC^*=(AA^*,BB^*)\in\Rad\mf{A}$ and $CC^*$ is
nilpotent element, so $AA^*$ and $BB^*$ are nilpotent and positive.
Therefore $C=(0,0)$, and the algebra $\mf{A}$ is semisimple.
On the other hand, the algebra $\mf{A}$ is local as an algebra
of endomorphisms of indecomposable representation.
Hence, $\mf{A}\cong \mb{C}$.
\end{remark}

\begin{lemma}\label{le01}
Let
$Z=\begin{bmatrix}z_{ij}\end{bmatrix}\substack{i=\ov{1,m}\\j=\ov{1,n}},\;
W=\begin{bmatrix}w_{ij}\end{bmatrix}\substack{i=\ov{1,m}\\j=\ov{1,n}}$
be matrices over the field $\mb{C}$, which has equal positive lengths
($|\oa{x}|,\, |y^\downarrow|$)
of corresponding rows and corresponding columns.
Let $A=\diag\{\seq{a}{m}\},\; B=\diag\{\seq{b}{n}\}$ be matrices over $\mb{R}$,
$a_i>0,\, b_j>0$ for $i=\ov{1,m},\, j=\ov{1,n}$, and let $AZ=WB$.
Then $Z=W$.
\end{lemma}

\emph{Proof.} Denote $K$~--- the number of nonzero elements in
the matrices $Z$ and $W$ ($AZ=WB$ implies that this number is the same
for $Z$ and $W$). Obviously, $K\geq \max(m,n)$, since the matrices $Z,\,W$
have no zero rows or columns.

We will use induction by the triples of numbers $(m,n,K)$:
assume that
$(m_1,n_1,K_1)<(m_2,n_2,K_2)$ if $m_1\leq m_2,\, n_1\leq n_2$ and
at least one inequality is strict, or if $m_1=m_2,\, n_1=n_2$ but $K_1<K_2$.

\medskip

a) The base of induction is obtained when $m=1$, either $n=1$, or $K=\max(m,n)$.

Let $m=1$. From
$[a_1z_{11},a_1z_{12},\ldots,a_1z_{1n}]=[w_{11}b_1,w_{12}b_2,\ldots,w_{1n}b_n]$ and
$|z_{1j}|=|w_{1j}|$ (the last follows from the equality of lengths of corresponding
columns of $Z$ and $W$), where $z_{1j}\neq 0$ and $w_{1j}\neq 0$
for $j=\ov{1,n}$, we obtain that $a_1=b_1=b_2=\cdots =b_n$, and then
$[z_{11},z_{12},\ldots,z_{1n}]=[w_{11},w_{12},\ldots,w_{1m}]$.

The case $n=1$ is analogous.

Now let for the sake of definiteness $1<m\leq n$. Then $K\geq n$.

Let $K=n$. In this case each column of the matrices $Z$ and $W$ contains
precisely one nonzero element.
Consider corresponding nonzero elemnts $z_{ij}$ and
$w_{ij}$ of the matrices $Z$ and $W$.
$AZ=WB$ implies $a_iz_{ij}=w_{ij}b_j$.
Since $|z_{ij}|=|w_{ij}|$ then $a_i=b_j$, and then $z_{ij}=w_{ij}$ and,
consequently, $Z=W$.

\medskip

b) Let $K>n\geq m>1$.

Since matrices $Z$ and $W$ has equal lengths of corresponding rows and
columns (by the assumption), and matrices $AZ$ and $WB$~---
because of the equality of matrices,
we obtain the following relations:
\begin{gather}\label{eq09}
a_i^2\left(|z_{i1}|^2+\cdots +|z_{in}|^2\right)=
|w_{i1}|^2b_1^2+\cdots +|w_{in}|^2b_n^2,\nonumber\\
|z_{i1}|^2+\cdots +|z_{in}|^2=
|w_{i1}|^2+\cdots +|w_{in}|^2,\nonumber\\
i=\ov{1,m};
\end{gather}
\begin{gather}\label{eq10}
a_1^2|z_{1j}|^2+\cdots +a_m^2|z_{mj}|^2=
\left(|w_{1j}|^2+\cdots +|w_{mj}|^2\right)b_j^2,\nonumber\\
|z_{1j}|^2+\cdots +|z_{mj}|^2=
|w_{1j}|^2+\cdots +|w_{mj}|^2,\nonumber\\
j=\ov{1,n};
\end{gather}
\begin{equation}\label{eq11}
a_i^2=\frac{|w_{i1}|^2}{\sum_{j=1}^n |w_{ij}|^2}b_1^2+\cdots
+\frac{|w_{in}|^2}{\sum_{j=1}^n |w_{ij}|^2}b_n^2,\; i=\ov{1,m};
\end{equation}
\begin{equation}\label{eq12}
b_j^2=a_1^2\frac{|z_{1j}|^2}{\sum_{i=1}^m |z_{ij}|^2}+\cdots
+a_m^2\frac{|z_{mj}|^2}{\sum_{i=1}^m |z_{ij}|^2},\; j=\ov{1,n}.
\end{equation}
In the right parts of equalities in \eqref{eq11}-\eqref{eq12} the
$a_i^2,\,b_j^2$ has coefficients that $\leq 1$.

We will assume (perhaps, after permutations of rows and columns) that
$a_1\leq a_2\leq \cdots \leq a_m$ and
$b_1\leq b_2\leq \cdots \leq b_n$.

Let $w_{1r}$ be the first and $w_{1s}$ be the last nonzero element in
the first row in the matrix $W$, $1\leq r\leq s\leq n$. Then
$$
a_1^2=\sum_{k=r}^s \frac{|w_{1k}|}{\sum_{j=r}^s |w_{1j}|^2}b_k^2\geq
b_r^2\sum_{k=r}^s \frac{|w_{1k}|^2}{\sum_{j=r}^s |w_{1j}|^2}=b_r^2.
$$
On the other hand,
$$
b_1^2=\sum_{l=1}^m \frac{|z_{l1}|}{\sum_{i=1}^m |z_{i1}|^2}a_i^2\geq
a_1^2\sum_{l=1}^m \frac{|z_{l1}|^2}{\sum_{i=1}^m |z_{i1}|^2}=a_1^2,
$$
\noindent i.~e. $a_1^2\geq b_r^2\geq b_1^2\geq a_1^2$.
Therefore $a_1^2=b_r^2$, and then $AZ=WB$ implies $z_{1r}=w_{1r}$.

Substitute elements $w_{1r}$ and $z_{1r}$ in $W$ and $Z$ by zeros, and
if after that matrices $W,\,Z$ contain zero rows or zero columns, eject them.
As a result we obtain matrices $\wh{W},\, \wh{Z}$ satisfying the conditions
of lemma and havind either the less number of zero elements
with the same dimension $(m\times n)$, or matrices with less
number of rows or columns (and at the same time with less number of
nonzero elements).
By the assumption of induction $\wh{W}=\wh{Z}$, and then, of course, $W=Z$.

\bigskip

Using the lemma, we will prove the following theorem:

\begin{theorem}\label{th01}
Let $Q$ be separated single quiver and $T$ be its indecomposable (in the
category $\ReposQH$) orthoscalar representation. Then $T$ is indecomposable
Schur representation in the category $\Rep Q$.
\end{theorem}

\emph{Proof.} Let $C=(A,B)$ be endomorphism of representation $T$ in the
category $\Rep Q$, i.~e. $AT=TB$. We may consider, that $A$ and $B$ are
invertible matrices, adding to $A$ and $B$, if it is necessary,
the apropriate scalar matrices with the same scalar on the diagonal.

Let $A=XU,\; B=VY$ be polar decompositions of the matrices $A$ and $B$,
where $U,\,V$ are unitary, $X,\,Y$ are positive nonsingular matrices
(here we may assume that $U,\,V,\,X,\,Y$ has the same block-diagonal structure,
as $A$ and $B$). Let besides
$$
X=U_1^*\wt{X}U_1,\; Y=V_1\wt{Y}V_1^*,
$$
\noindent where $U_1,\,V_1$ are unitary matrices, $\wt{X},\,\wt{Y}$ are
diagonal matrices with positive numbers on the diagonal. Then
$$
U_1^*\wt{X}U_1UT=TVV_1\wt{Y}V_1^*
$$
\noindent or
$$
\wt{X}(U_1UTV_1)=(U_1TVV_1)\wt{Y}.
$$
Since the lengths of corresponding rows and columns of matrices $U_1UTV_1$
and $U_1TVV_1$ are equal because of orthoscalarity of representation $T$,
we obtain by the lemma~\ref{le01}
$$
U_1UTV_1=U_1TVV_1,
$$
\noindent and then, after reduction,
$$
UT=TV.
$$
Since $U,\,V$ are unitary matrices, then by remark~\ref{re02} $(U,V)$ is an
endomorphism of indecomposable in $\ReposQH$ orthoscalar representation, and then
by remark~\ref{re04} the matrices $U$ and $V$ are scalar with the same scalar
on the diagonal.

Reduce the equality
$$
XUT=TVY
$$
\noindent by the scalar matrices $U$ and $V$, obtain
$$
XT=TY,
$$
\noindent where matrices $X$ and $Y$ are self-adjoint.

Hence by remark~\ref{re03} $(X,Y)$ is endomorphism of orthoscalar indecomposable
representation $T$ in the category $\ReposQH$ and (by remark~\ref{re04})
$X$ and $Y$ are scalar with the same scalar on the diagonal.
As a result the same will also be products of scalar matrices $XU=A$ and $VY=B$,
and so the representation $T$ is Schur representation in the category $\Rep Q$.

Since, obviously, decomposable representation has nonscalar
endomorphism, the representation $T$ is indecomposable in $\Rep Q$.

\bigskip

\begin{remark}\label{re05}
The condition that a quiver is single is inessential, the lemma~\ref{le01}
is easy generalized on the case of a quiver with multiple arrows.
\end{remark}

\begin{remark}\label{re06}
If a quiver $Q$ is not separated, the theorem~\ref{th01} will not be
true. Indeed, let quiver $Q$ be the loop

\bigskip
\centerline{\includegraphics[scale=0.4]{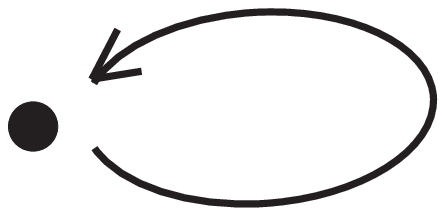}}
\bigskip

Let $T$ be its indecomposable orthoscalar representation with matrix
$$
T=\begin{bmatrix}
1&1\\
0&-1
\end{bmatrix}.
$$
In this case the condition orthoscalarity is \cite{KruRoi1}
$$
TT^*+T^*T=3I.
$$
The indecomposability of $T$ in $\ReposQH$ is checked as follows:
it easy to calculate that $CT=TC$ and $CT^*=T^*C$ imply,
that $C$ is a scalar matrix.

Let $A=\begin{bmatrix}3&1\\ 0&1\end{bmatrix}$; it is easy to check that
$AT=TA$, so that $A$ defines nonscalar invertible
endomorphism of representation $T$ in $\Rep Q$
(and so representation $T$ is decomposable in $\Rep Q$).
\end{remark}

\bigskip

{\bf 3.} Set $\mc{P}_n$ be the following $*$-algebra:
$$
\mc{P}_n=\mb{C}\langle
\seq{p}{n}\,|\, p_i=p_i^2=p_i^*
\rangle,
$$
\noindent and $\mc{K}(\mc{P}_n)$ be the category, which objects are
$*$-representations $\pi$ of the
algebra $\mc{P}_n$ in the category of finite-dimensional
Hilbert spaces $\mc{H}$
($\pi(p_i)=P_i,\; P_i: H_0\to H_0,\; P_i^2=P_i^*=P_i$), and
a morphism $C_0$ from representation $\pi$ in space $H_0$
to representation $\wt{\pi}$ in space $\wt{H}_0$
is defined as a linear map
$C_0: H_0\to \wt{H}_0$ with the property
\begin{equation}\label{eq13}
C_0P_i=\wt{P}_iC_0P_i,\; i=\ov{1,n}.
\end{equation}
One can regard the category $\mc{K}(\mc{P}_n)$ \cite{EnoWat1} as the category
of collections of $n$ subspaces of finite-dimensional Hilbert spaces.
If $\pi$ is a representation of the algebra
$\mc{P}_n$ in $H_0$ ($\pi(p_i)=P_i$), then it associates with $n$
subspaces $(\seq{H}{n})$ of space $H_0$;
here $H_i=\Imm P_i,\; i=\ov{1,n}$. A morphism from
$(\seq{H}{n})$ to $(\seq{\wt{H}}{n})$ is defined naturally
as linear map $C: H_0\to \wt{H}_0$, for which
$C(H_i)\subset \wt{H_i}$ for $i=\ov{1,n}$. The category
$\mc{K}(\mc{P}_n)$ and the category of collections of $n$ subspaces
of finite-dimensional Hilbert spaces are equivalent.
Papers \cite{GelPon1,Bre1,NazRoi1} and others are dedicated to the classification
of collections of $n$ subspaces in linear spaces.

Representation $\pi$ algebra $\mc{P}_n$ is called orthoscalar, if there exist
such positive numbers $\seq{\alpha}{n}$ that
$$
\sum_{i=1}^n \alpha_i\pi(p_i)=I_0
$$
\noindent ($I_0$ is identity operator in space $H_0$).
Denote as $\mc{L}(\mc{P}_n)$ the full subcategory of orthoscalar representations
in $\mc{K}(\mc{P}_n)$. Denote as $Q_n$ the quiver

\bigskip
\centerline{\includegraphics[scale=0.4]{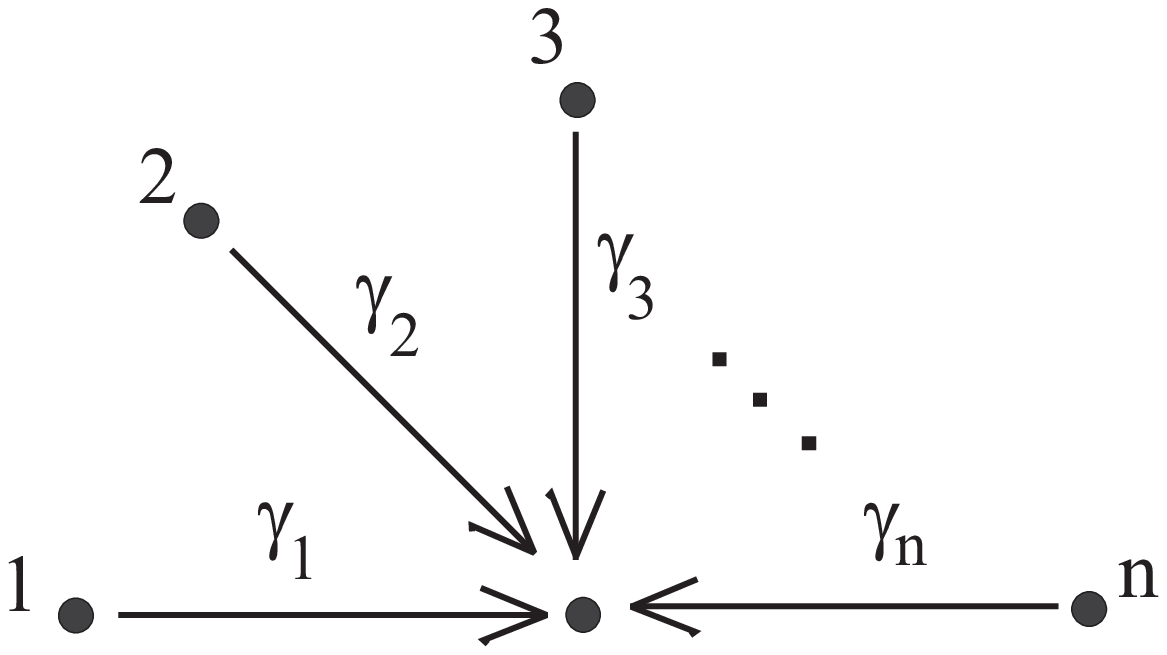}}
\bigskip

Let $\mc{K}(Q_n)$ be the full subcategory in $\Rep Q_n$ of those representations
$T$ for which $T(\gamma_i)^*T(\gamma_i)$
is a scalar operator, different from zero in space $T(i)$, and
$\ReposQnH$ be the full subcategory of orthoscalar representations in
the category $\mc{K}(Q_n)$.

\begin{lemma}\label{le02}
Categories $\mc{K}(Q_n)$ and $\mc{K}(\mc{P}_n)$ are equivalent.
Categories $\ReposQnH$ and $\mc{L}(\mc{P}_n)$ are equivalent.
\end{lemma}

\emph{Proof.} Construct a functor $F: \mc{K}(Q_n)\to \mc{K}(\mc{P}_n)$.
Let $T$ be a representation from $\mc{K}(Q_n)$ and
$T(\gamma_i)^*T(\gamma_i)=\alpha_iI_i$. Put
$\frac{1}{\sqrt{\alpha_i}}T(\gamma_i)=\Gamma_i$, then $\Gamma_i^*\Gamma_i=I_i$ and
$\Gamma_i\Gamma_i^*=P_i$ is orthogonal projection.
Let representation $\pi$ of the algebra $\mc{P}_n$
is defined by equalities $\pi(p_i)=P_i$. Put $F(T)=\pi$.

Let $C=\{C_i\}_{i=\ov{0,n}}: T\to\wt{T}$ be a morphism in the category $\mc{K}(Q_n)$:
$$
C_0T(\gamma_i)=\wt{T}(\gamma_i)C_i,\; i=\ov{1,n}.
$$
Then
\begin{equation}\label{eq14}
\sqrt{\alpha_i}C_0\Gamma_i=\sqrt{\wt{\alpha}_i}\wt{\Gamma}_iC_i,
\end{equation}
\noindent and since $\wt{\Gamma}_i^*\wt{\Gamma}_i=\wt{I}_i$ then
\begin{equation}\label{eq15}
C_i=\sqrt{\frac{\alpha_i}{\wt{\alpha}_i}}\,\wt{\Gamma}_i^*C_0\Gamma_i.
\end{equation}
After substitution of $C_i$ in \eqref{eq14} we obtain
$$
C_0\Gamma_i=\wt{\Gamma}_i\wt{\Gamma}_i^*C_0\Gamma_i,
$$
\noindent and, after multiplying of equality on the right by $\Gamma_i^*$,
$$
C_0P_i=\wt{P}_iC_0P_i,
$$
\noindent i.~e. $C_0: \pi\to\wt{\pi}$ is a morphism in the
category $\mc{K}(\mc{P}_n)$.

Put $F(C)=C_0$, then $F$ is a functor from $\mc{K}(Q_n)$
to $\mc{K}(\mc{P}_n)$.

Let $C_0: \pi\to\wt{\pi}$ be a morphism in $\mc{K}(\mc{P}_n)$ where
$C_0$ is a linear map from $H_0$ to $\wt{H}_0$. Construct maps
$C_i: H_i\to \wt{H}_i$ where $H_i=\Imm P_i,\; \wt{H}_i=\Imm \wt{P}_i$
by the formulas \eqref{eq15}. If $\Gamma_i,\,\wt{\Gamma}_i$ are
natural embeddings of $H_i$, resp. $\wt{H}_i$ in space $H_0$,
resp. $\wt{H}_0$
($P_i=\Gamma_i\Gamma_i^*,\; \wt{P}_i=\wt{\Gamma}_i\wt{\Gamma}_i^*$),
then
$$
\wt{T}(\gamma_i)C_i=
\sqrt{\wt{\alpha}_i}\wt{\Gamma}_iC_i=
\sqrt{\wt{\alpha}_i}\wt{\Gamma}_i\sqrt{\frac{\alpha_i}{\wt{\alpha_i}}}
\wt{\Gamma}_i^*C_0\Gamma_i=
\sqrt{\alpha_i}\wt{P}_iC_0P_i\Gamma_i=
\sqrt{\alpha_i}C_0P_i\Gamma_i=
\sqrt{\alpha_i}C_0\Gamma_i=
C_0T(\gamma_i),
$$
\noindent i.~e. $C=\{C_i\}_{i=\ov{0,n}}$ is such morphism in $\mc{K}(Q_n)$ that
$F(C)=C_0$. Since, using formulas \eqref{eq15} $C_i$ are defined uniquely by $C_0$,
$F$ is a complete and univalent functor. Each
representation in $\mc{K}(\mc{P}_n)$ is obtained from a certain representation
in $\mc{K}(Q_n)$, therefore functor $F$ is an equivalence of categories.

Resrtiction of $F$ on $\ReposQnH$ gives equivalnce of $\ReposQnH$ and
$\mc{L}(\mc{P}_n)$.

The following therem is a corollary from
lemma~\ref{le02} and theorem~\ref{th01}:

\begin{theorem}\label{th02}
If $(\seq{H}{n})$ is a collection $n$ subspaces of finite-dimensional Hilbert
space, and corresponding collection of $n$ orthogonal projections
$(\seq{P}{n})$ define an indecomposable representation in the category
of orthoscalar representations $\mc{L}(\mc{P}_n)$
($\sum\limits_{i=1}^n \alpha_iP_i=I_0$ for a certain collcetion of
numbers $\alpha_i>0,\; i=\ov{1,n}$), then this collection of orthogonal projections
is a Schur object in the category $\mc{K}(\mc{P}_n)$, and a corresponding
collection of subspaces is an indecomposable Schur object in the
category of collections of $n$ subspaces of finite-dimensional Hilbert space.
\end{theorem}

\end{document}